%%%%%%%%%%%%%%%%%%%%%%%%%%%%%%%%%%%%%%%%%
%
%   Lengths of Periods and Seshadri
%   Constants of Abelian Varieties
%
%         Robert Lazarsfeld
%
%              AMS-TeX
%
%
%       version of June 13, 1996
%
%
%%%%%%%%%%%% RKL AMS-TeX SETTINGS %%%%%
\input amstex
\documentstyle{amsppt}
\magnification \magstep1
\parskip 11pt
\parindent .3in 
\pagewidth{5.25in} 
\pageheight{7.4 in}

%%%%%%%%%%%%%%%%%%%%%%%%%%%%%%%%%%%%%%%%%
%
%%%%  COMMON  RKL  DEFINITIONS
%
\def \bl{\vskip 11pt} 

\def \ni{\noindent}
\def \P{\bold{P}}

\def \R{\bold{R}}
\def \C{\bold{C}} 
\def \O{\Cal{O}}
\def \eps{\epsilon}
\def \wstd {\omega_{\text{std}}}
\def \l {\lambda}
\def \tbar {\overline{\tau}}
\def \w {\omega}
\def \s {\sigma}
\def \wbar{\overline \w_L}
\def \lra{\longrightarrow}
%
%%%%%%%%%%%%%%%%%%%%%%%%%%%%%%%%%%%  END OF STANDARD
%                                    RKL DEFINITIONS %%%%%%%%%

\centerline{\bf LENGTHS OF PERIODS AND SESHADRI CONSTANTS}
\vskip 2pt
\centerline{\bf OF ABELIAN VARIETIES}

\vskip 15pt
\centerline{Robert LAZARSFELD\footnote{Partially supported by N.S.F.
Grant DMS 94-00815}}
\bl

\ni {\bf Introduction}

The purpose of this note is to point out an elementary but  somewhat
surprising connection between the work of Buser and Sarnak \cite{BS} on
lengths of periods of abelian varieties and the  Seshadri constants
measuring the local positivity of theta divisors. The link
is established via symplectic blowing up, in the spirit of \cite{McDP}.

We start by recalling the definition of Seshadri constants. Let $X$ be a
smooth complex projective variety, let
$L$ be an ample line bundle on $X$, and fix a point $x \in X$. Consider
  the blowing-up
$$f : Y = \text{Bl}_x(X) \lra X$$
of $X$ at $x$, with exceptional divisor $E = f^{-1}(x) \subset Y$. Then
for $0 < \eps \ll 1$ the cohomology class $f^* c_1(L) - \eps \cdot [E]$ will
lie in the K\"ahler cone of $Y$. As a measure of how positive $L$ is
locally near $x$ we ask in effect how large we can take $\eps$ to be while
keeping the class in question positive. More precisely, set
$$\eps(L,x) = \sup \left \{  \ \eps \ge 0 \mid f^*c_1(L)   - \eps \cdot [E] \
\text{ is nef} \ \right
\}.$$ Here $f^* c_1(L) - \eps \cdot [E]$ is considered as an $\R$-divisor
class on
$Y$, and to say that it is nef means that $\int_{C^\prime} f^* c_1(L) \ge
\eps(  E
\cdot C^\prime) $ for every irreducible algebraic curve $C^\prime
\subset Y$. \footnote{Recall that a theorem of Kleiman characterizes  the
nef cone as the closure of the ample cone.}  We refer to
\cite{Dem, \S6} or \cite{EKL, \S1} for further discussion and alternative
characterizations. Introduced by Demailly in \cite{Dem}, these
Seshadri constants have attracted considerable interest in recent
years. The main result of \cite{EKL} states that if
$X$ has dimension $n$, then at a very general point $x \in X$ one has the
universal lower bound $\eps(L, x)
\ge \frac{1}{n}$ (cf. also \cite{KS}). Some more refined results when
$X$ is a surface appear in \cite{EL},  \cite{S} and \cite{Xu2}, but 
except in the simplest examples Seshadri constants have
proven very difficult to control with any precision. 
We propose here to study these invariants when the
ambient manifold  is an abelian variety.

Suppose then that $(A, \Theta)$ is a principally polarized abelian
variety of dimension $g$, i.e. that $A$ is a complex torus, and
that $\Theta
\subset A$ is an ample divisor with $h^0(A, \O_A(\Theta)) = 1$.  
Since $A$ is homogeneous, the Seshadri constants
$\eps(\O_A(\Theta),x)$ are
independent of
$x \in A$, and we denote  their common value by $\eps(A, \Theta)$ or
simply $\eps(A)$.
 One has the elementary upper bound
$$\eps(A) \le \root \uproot 2 g \of { g!} $$ (cf. \cite{EKL,
1.8}). Nakamaye
\cite{N} has shown that $\eps(A, \Theta) \ge 1$, with equality iff $(A,
\Theta)$ is the product of an elliptic curve and an abelian variety of
dimension $g-1$.

Our goal is to relate the Seshadri constant $\eps(A)$ to a
metric invariant of $(A, \Theta)$. As usual, write
$A$ as a quotient
$$ A = V / \Lambda$$ of its universal covering, so that $V \cong \C^g$,
and $\Lambda \subset V$ is a lattice in $V$. The
principal polarization
$\Theta$ determines a positive definite Hermitian form $H$ on $V$ (cf.
\cite{LB,Chapter 2}),  and following \cite{BS} we
define 
$$m(A) = m(A, \Theta) = \min_{x\in \Lambda - \{ 0 \} } H(x,
x).$$ Thus $m(A)$ is the square of the minimal length (with respect
to  $H$) of a non-zero lattice  vector. This is the analogoue for
abelian variety period lattices of  an invariant
 familiar in connection with sphere packings and the geometry of
numbers (cf.
\cite{O}).  Buser
and Sarnak  study the maximum value of $m(A)$ as $A$
varies over the moduli space $\Cal A_g$ of principally polarized abelian
varieties, and they show (\cite{BS, \S2}) that there exist p.p.a.v.'s
$(A,
\Theta)$ for which
$$ m(A)  \ge \frac{1}{\pi} \left ( 2 g ! \right
)^{1/g} . \tag BS1  $$
The most surprising result of \cite{BS} is that if $C$ is a smooth
projective algebraic curve of genus $g \ge 2$, and $(J(C), \Theta_C)$ is
its polarized Jacobian, then one has the upper bound
$$m(J(C))
\le
\frac{3}{\pi}\log (4g + 3). \tag BS2$$
In other words, for $g \gg 0$ a Jacobian has a period of  unusually short
length.

Our main result states that the Seshadri constant of $A$ is
bounded below in terms of the mimimal length of a period:
\proclaim {\bf Theorem} One has the inequality
$$\eps(A) \ge \frac{\pi}{4} m(A). $$
\endproclaim
\ni Note that in general (and maybe always) the inequality is strict,
as one sees already in the one dimensional case. 

	This inequality has a number of pleasant consequences. In the
first place, combining the Theorem with the bound (BS1) of Buser and
Sarnak, we obtain the
\proclaim{Corollary} Let $(A, \Theta)$ be a very general principally
polarized abelian variety.  Then
$$\eps(A) \  \ge  \ \frac{ 2^{\frac{1}{g}}}{4} \root \uproot 4 g \of
{g!}  
 \  \  \approx   \ \ \frac{g}{4e}.$$
\endproclaim
\ni The hypothesis on $A$ means that the inequality
is valid off the union of countably many proper subvarieties of the
moduli space $\Cal A_g$. In the approximation, which holds for
$g
\gg 0$, we are ignoring the factor of
$2^{1/g}$. Observe that this lower bound differs from
the upper bound
$\eps(A) \le \root \uproot 2 g \of { g!} $ by  a factor of less than $4$.
It would be interesting to know whether  $\eps(A_{\text{very general}}) =
(g!)^{1/g}$ for large $g$.\footnote{Proposition 3 of \cite{S}
asserts that equality never holds, but the proof is erroneous (a
circumstance for which the present author must share some
culpability).}

Now let $C$ be a compact Riemann surface 
of genus $g \ge 2$, and as above let
$(J(C),
\Theta_C)$ be its polarized Jacobian. It is rather easy to
obtain upper bounds on the Seshadri constants of $J(C)$:
\proclaim {\bf Proposition} {\rm (i)}. One has
$$\eps(J(C), \Theta_C) \le \sqrt  g.$$

\ni {\rm (ii).} Suppose that $C$ 
can be expressed as  
a $d$-sheeted branched covering $\phi : C \lra \P^1$. Then 
$$\eps(J(C), \Theta_C) \le \frac{gd}{g + d - 1}.$$
\endproclaim
\ni For hyperelliptic curves (when $d = 2$), the inequality (ii) was
established by Steffens 
\cite{S}. Combining  statement (i) and the Theorem, one arrives at an 
elementary new
 proof that  a Jacobian has a period of  small length,
although the specific inequality that comes out is not as strong as (BS2).
On the other hand, we see from (ii) that if $C$ is a $d$-sheeted covering
of $\P^1$,
then in fact
$$m(J(C)) \le \frac{4d}{\pi}.$$
This seems to be new. 
In the other direction, Buser and Sarnak construct examples of curves to
show that the supremum of $m(J(C))$ on the moduli space $\Cal{M}_g$
is $\ge c \cdot \log(g)$, where $c$ is a small positive constant.
Hence   the Seshadri constant $\eps(J(C_{\text
{very general}}))$ of the Jacobian of a very general curve
satisfies the same inequality (with a slightly
different constant). It would be interesting to know how
$\eps(J(C_{\text
{very general}}))$
actually grows with $g$. It is also tempting to wonder to what extent
  small
Seshadri constants might characterize Jacobians   among
all irreducible p.p.a.v.'s.

I am grateful to L. Ein, 
M. Nakamaye and M. Thaddeus for valuable discussions.
I'd also like to acknowledge my debt to
 the papers \cite{Deb} and \cite{Xu1},  
through which I became  aware of \cite{BS} and \cite{McDP} respectively. 

\vskip 10pt
\ni{\bf \S 1.  Local Positivity on Abelian Varieties}
\vskip 5pt
The theorem is a simple consequence of the construction of the
symplectic blowing up of a point, as explained for example in the
paper \cite{McDP} of McDuff and Polterovich. The basic point, which is
implicit in \cite{McDP}, is a relation between Seshadri constants and
radii of symplectically embedded holomorphic balls. This connection was
exploited in a related but more sophisticated manner in \cite{McDP}.

We start by fixing notation. In  
$\C^n$ with coordinates
$z_j =  x_j + i   y_j$, denote by
\def \wstd {\omega_{\text{std}}}
$$\wstd = \sum dx_j \wedge dy _j = \frac{i}{2} \sum dz_j \wedge  d
\overline z_j$$ 
 the standard symplectic form. Write  $B(\l) \subset
\C^n$ for the open ball of radius $\l$ centered at the origin:
$$B(\l) = \{ z \in \C^n \mid |z|^2 < \l^2 \}.$$ We view
$B(\l)$ as a complex manifold, and also as a symplectic manifold via
$\wstd$. 

Now let $X$ be a smooth projective variety of dimension $n$, $L$ an ample
line bundle on
$X$, and 
$\w_L$ a K\"ahler form\footnote{I.e. a closed positive $(1,1)$ form.} on
$X$ representing $c_1(L)$. We view $(X, \w_L)$ as a symplectic
manifold. Given $x \in X$ we define a real number $\l(x) = \l(\w_L,x)
\ge 0$ by looking for the largest radius $\l>0$ for which there exists a
holomorphic and symplectic embedding
$$j = j_{\l}  : (B(\l) , \wstd)  \hookrightarrow (X , \w_L) \ \
\text{with}
\ \ 0 \mapsto x. \tag *$$
More precisely, if there is no $\l > 0$ for which an embedding (*)
exists, set $\l(x) = 0$. Otherwise, put
$$\l( \omega_L, x) = \sup \left \{ \l > 0 \mid \ \exists  \
\text{holomorphic and symplectic }  j_\l \text{ as in (*) } \right \}.$$
\proclaim{Main Lemma}  One has
the inequality
$$\eps(L, x) \ge \pi \l( \w_L, x)^2.$$
\endproclaim
\ni By way of proof, it would probably be almost enough just to refer to
\cite{McDP}, (5.1)-(5.3). But since the lemma isn't stated
there explicitly, and since it involves some ideas that are not
 standard algebro-geometrically, we will summarize the argument
for the benefit of the reader in \S 2. In the meantime, we  grant the
lemma. 

The rest of the proofs
are quite immediate:

\demo{Proof of Theorem} Let $\pi : V \lra A$ be the universal covering,
and as above let $H$ be the Hermitian form on $V$ determined by
$\Theta$.
In the natural way, we may view the imaginary part $\w = \text{im}  \
H$ as a symplectic form on $V$, which is in fact the pull-back
 $\w = \pi^*
\w_{\Theta}$ of a K\"ahler form $\w_{\Theta}$ on $A$ representing
$c_1(\O_A(\Theta))$. We   fix a basis of $V$ with respect to which 
$H$ is the standard Hermitian form
$$H(v, w) \ =  \ ^t \overline v \cdot \ w$$
on $\C^g$. Then taking $z_j$ to be the corresponding complex
coordinates, one has 
$$ \w = \pi^* \w_{\Theta} = \wstd, \tag *$$ and $H(x,x) = | x |^2$ is
just the usual Euclidean length. In particular,
$$m(A) = \min_{x \in \Lambda - \{
0 \}} \left \{   \ |x|^2  \  \right \}.$$

Now let $\l = \sqrt {m(A)}/ 2$. 
Then given any two points $x , y \in
B(\l)$ one has $|x - y| <  2 \l = \sqrt{m(A)}$. Therefore no two points
of
$B(\l)$ are congruent (modulo $\Lambda$), and consequently the
composition
$$j_\l : B(\l) \hookrightarrow V \overset \pi \to \lra A$$ is an
embedding. But $j_\l$ is of course holomorphic, and thanks to (*) it is
symplectic as well. Therefore
$$\l( \w_\Theta , 0) \ge  \frac{\sqrt{m(A)}}{2},$$
and the Theorem follows from the Main  Lemma. \qed \enddemo

\demo{Proof of Proposition}  We assume to begin with that $C$ is
non-hyperelliptic. Consider the subtraction map
$$s : C \times C \lra J(C), \ \ \ (x, y) \mapsto \O_C(x-y) \in 
Pic^0(C) = J(C),$$ and let $\Sigma \subset J(C)$ be its image. It is
elementary and well known (cf. \cite{ACGH, pp. 223, 263}) that if $C$ is
non-hyperelliptic, then
$s$ is an isomorphism off the diagonal $\Delta \subset C \times C$,
and blows $\Delta$ down to the origin $0 \in\Sigma$, which is a point
of multiplicity $2g - 2$. Moreover $\Delta$ is the scheme-theoretic
inverse image of  the singular point $0 \in \Sigma$. The required
inequalities will follow from some computations in the intersection ring of
$C \times C$. To this end, let  $F_1, F_2\subset C \times C$ 
be the preimages of
a point of $C$ under the two projections. Then working with numerical
equivalence of divisors,  one checks that
$$s^*(\Theta) \equiv (g-1)(F_1 + F_2) + \Delta$$
(cf. \cite{R}). It follows with a calculation that the degree of
$\Sigma$ with respect to $\Theta$ is
$$\aligned
\text{deg}_{\Theta}(\Sigma) = \Theta^2 \cdot \Sigma &= \big (
(g-1)(F_1 + F_2) + \Delta \big ) ^2 \\
&= 2g(	g-1).
\endaligned$$
Then by \cite{Dem, (6.7)}:
$$
\aligned
\eps(J(C), \Theta) &\le \sqrt{ \frac{\text{deg}_\Theta(\Sigma)
}{\text{mult}_0\Sigma}} \\
&= \sqrt{ \frac{2g(g-1)}{2g - 2}} \\
&= \sqrt{g}.
\endaligned $$

Turning to statement (ii), let $L = \phi^* \O_{\P^1}(1)$.
 Then there is an effective divisor $\Gamma
\subset C \times C$ with $$\Gamma \in |pr_1^* L \otimes
pr_2^*L \otimes \O_{C\times C}(-\Delta)|.$$
Geometrically, for instance, we may realize $\Gamma$ as the closure of 
$\Gamma_0 = \{ (x, y ) \mid x\ne y , \phi(x) = \phi(y) \}$. Now if
$\eps = \eps(J(C), \O_J(\Theta))$, then $s^*(\Theta) - \eps \cdot
\Delta$ is nef on $C \times C$. Hence
$$\aligned
\Gamma \cdot \left ( s^*(\Theta) - \eps \Delta \right ) 
&=  \big ( d(F_1 + F_2) - \Delta  \big ) \cdot 
 \big ( (g-1)(F_1 + F_2) + (1 - \eps) \Delta  \big ) \\ 
&\ge 0, 
\endaligned$$
and with another calculation this leads to the second assertion of the
Proposition. Finally, if $C$ is hyperelliptic the only thing that
needs proof is statement (ii) with $d = 2$, and this follows by
looking at the image in $J(C)$ of the curve $\Gamma$ just constructed.
\qed

\vskip 15pt

\ni{\bf \S2. Sketch of Proof of Main Lemma}

Finally, for the benefit of  readers
not versed in symplectic matters, we outline  the proof of the Main
Lemma. We follow
\cite{McDP},
\S5  (also pp. 414 ff), quite closely. The essential point, which seems
to go back at least as far as  \cite{GS}, is to
construct explicitly a   K\"ahler form on the blow-up Bl$_0(\C^n)$ 
which agrees with the standard form off a ball of specified radius. The
presence of a symplectically embedded holomorphic ball allows one to
carry over the local construction to a global setting, and then the
inequality of the Main Lemma  follows from the positivity of the form so
constructed.

Turning to the details, let
$$V \subset \C^n \times \P^{n-1}$$
be the blowing up of $0 \in \C^n$, embedded in the usual way as an
incidence correspondence. Write
$$f: V \lra \C^n , \ \ \ \ q : V \lra \P^{n-1}$$
for the projections, so that $f$ is the blowing-up, and $q$ realizes
$V$ as the total space of the line bundle $\O_{\P^{n-1}}(-1)$. Denote
by $V(\l)$ the inverse image of the ball $B(\l) \subset \C^n$:
$$V(\l) = f^{-1} B(\l) \subset V,$$ 
so that $V(\l)$ is an open neighborhood of the exceptional divisor
$E = \P^{n-1} \subset V$. Finally, let 
$\s$ be the usual Fubini-Study K\"ahler form on $\P^{n-1}$, normalized
so that $\int_{\P^1} \s  = \pi$, the integral being taken over a line
in $\P^{n-1}$. This normalization is chosen so that if $S = S^{2n-1}
\subset \C^n$ is the unit sphere, and $\kappa : S \lra \P^{n-1}$ is the
Hopf map, then $\kappa^* \s = \wstd | S$. 

The crucial ingredient is the following 
\proclaim {\bf Basic Local Construction} Fix $\l > 0$. Given any
small $\eta > 0$, there exists some $0 < \delta \ll 1$, plus a
K\"ahler form $ \ \tbar = \tbar(\l, \eta)$ on $V$ such that

{\rm (i)}. $\tbar = f^*(\wstd)  \ \text{ on }   \ V - \overline{
V\left (\l(1+\eta)\right )}$;

{\rm(ii)}. $\tbar = f^*(\wstd) + \l ^2 q^*(\s) \ \ \text{ on } \
V(\delta).$
\endproclaim
\ni In other words, $\tbar$ coincides with the standard K\"ahler form
on $\C^n$ off a ball of radius (a tiny bit larger than) $\l$, whereas
in a neighborhood of the exceptional divisor, we are ``twisting" by a
form representing $\pi \l^2 q^*  c_1(\O_{\P^{n-1}}(1))$. This is an
extremely  slight variant of \cite{McDP, (5.1)}, proved exactly as in 
\cite{McDP, (5.2), (5.3)}, and we refer the reader to the very clear
exposition there.\footnote{In brief, choose a monotone increasing
smooth function $\phi(r)$ such that $\phi(r) = \sqrt{\l^2 + r^2}$ for $0 <
r < \delta \ll 1$, and such that $\phi(r) = r$ for $r > \lambda(1 +
\eta)$, and then consider the smooth mapping 
$$F : \C^n - \{ 0 \} \lra \C^n, \ \ F(z) =  \frac{\phi(|z|)}{|z|} \cdot
z.$$ Then $\tbar = f^* F^* \wstd$, extended over
$E$ by (ii), has the required  properties.}  See also
\cite{McDS,
\S 6.2}.

Given this local construction, the proof of the main lemma is
rather evident. Let
$f : Y = \text{Bl}_x(X) \lra X$ be the blowing up of
$X$, with exceptional divisor $E \subset Y$, and fix any 
$\l < \l(\w_L,x)$. It is enough to show that 
$$\text{the  $\R$-divisor class  $f^*(c_1(L)) -
\pi \l^2 [E]$ is nef on $Y$.} \tag *
$$
To this end, fix
$0 < \eta \ll 1$ so that
$\l\cdot  (1 + 3 \eta) < \l(\w_L, x)$. We have a holomorphic and
symplectic embedding 
$$B(\l\cdot (1 + 3 \eta)) \hookrightarrow X, \tag **$$
and so for $\nu<  \l\cdot (1 + 3 \eta)$ we can view the local model
$V(\nu)$ as being embedded in $Y$ as a neighborhood of the exceptional
divisor. Thanks to property (i) and the fact that the embedding (**)
is symplectic,  the basic local construction guarantees   the existence
of a K\"ahler form $\wbar$ on
$Y$, agreeing with $\w_L$ off $V(\l(1 + 2 \eta))$, and being given by 
(ii) in a neighborhood $V(\delta)$ of $E$. Since $\wbar$ is K\"ahler,
and in particular positive, (*) will follow once we know that its
cohomology class satisfies
$$[\wbar] = f^* [\w_L] - \pi \l^2 [E] = f^* c_1(L) - \pi \l^2 [E].
\tag ***$$ But $\wbar - f^* \w_L$ is supported in a small neigborhood of
$E$, and then (***) follows easily   using (ii) and the
normalization of
$\s$. \qed

\vskip 10pt

\ni 
\vskip 5pt
\bl
\settabs\+University of Illinois at Chicago and stil extra 
\cr 
 
\+ Department of Mathematics \cr
\+ University of California,
Los Angeles \cr
\+ Los Angeles, CA  90024 \cr
\+ e-mail: rkl\@math.ucla.edu \cr

\end